\documentclass[12pt]{article}
\usepackage[latin9]{inputenc}
\usepackage{geometry}
\usepackage{mathrsfs}
\usepackage{amstext}
\usepackage{amssymb}

\makeatletter
 
\usepackage{amsfonts}
\usepackage{latexsym}
\usepackage{mathrsfs}

\title{Convexity and Osculation in Normed Spaces}
\author{Mark Mandelkern} 

\makeatother

\begin{document}
\newtheorem{defn}{Definition}[section] \newtheorem{prop}[defn]{Proposition}
\newtheorem{lm}[defn]{Lemma} \newtheorem{thm}[defn]{Theorem} \newtheorem{cor}[defn]{Corollary}
\newtheorem{ex}[defn]{Example} \setlength{\parindent}{15pt}
\begin{center}
\begin{Large} Convexity and Osculation in Normed Spaces \end{Large}\\
 \vspace{6mm}
 Mark Mandelkern 
\par\end{center}

\vspace{3mm}
 \pagenumbering{arabic} \setcounter{page}{1}
\begin{quotation}
\noindent \textbf{Abstract.} Constructive properties of uniform convexity,
strict convexity, near convexity, and metric convexity in real normed
linear spaces are considered. Examples show that certain classical
theorems, such as the existence of points of osculation, are constructively
invalid. The methods used are in accord with principles introduced
by Errett Bishop.\\
\end{quotation}
{\textbf{Introduction}} \label{sec:intro}

Contributions to the constructive study of convexity are found in
{[}B67{]}, {[}M83{]}, {[}S03{]} and, most extensively, in a recent
paper of Fred Richman {[}R07{]}. The present paper consists mainly
of comments related to the latter paper, and Brouwerian counterexamples
concerning various properties of convexity and the existence of certain
points.

\textit{Constructive mathematics.} A characteristic feature of the
constructivist program is meticulous use of the conjunction ``or".
To prove $``P\textnormal{ or }Q"$ constructively, it is required
that either we prove $P$, or we prove $Q$; it is not sufficient
to prove the contrapositive $\neg(\neg P\textnormal{ and }\neg Q)$.

To clarify the methods used here, we give examples of familiar properties
of the reals which are constructively \textit{invalid}, and also properties
which are constructively \textit{valid}. The following classical properties
of a real number $\alpha$ are constructively invalid: ``Either $\alpha<0$
or $\alpha=0$ or $\alpha>0"$, and ``If $\neg(\alpha\leq0\ ),$
then $\alpha>0"$. The relation $\alpha>0$ is given a strict constructive
definition, with far-reaching significance. Then, the relation $\alpha\leq0$
is defined as $\neg(\alpha>0)$. A constructively valid property of
the reals is the \textit{Constructive Dichotomy Lemma}: If $\alpha<\beta$,
then for any real number $\gamma$, either $\gamma>\alpha$, or $\gamma<\beta$.
This lemma is applied ubiquitously, as a constructive substitute for
the constructively invalid \textit{Trichotomy.} For more details,
see {[}B67{]}.

Applying constructivist principles when reworking classical mathematics
can have interesting and surprising results.\footnote{For an exposition of the constructivist program, see Bishop's Chapter
1, ``A Constructivist Manifesto", in {[}B67{]}; see also {[}M85{]}
and {[}R99{]}.}

\textit{Brouwerian counterexamples.} To determine the specific nonconstructivities
in a classical theory, and thereby to indicate feasible directions
for constructive work, \textit{Brouwerian counterexamples} are used,
in conjunction with \textit{omniscience principles.} A Brouwerian
counterexample is a proof that a given statement implies an omniscience
principle. In turn, an omniscience principle would imply solutions
or significant information for a large number of well-known unsolved
problems.\footnote{This method was introduced by L. E. J. Brouwer in 1908 to demonstrate
that use of the \textit{law of excluded middle} inhibits mathematics
from attaining its full significance.} Omniscience principles may be stated in terms of binary sequences
or properties of real numbers; we will have need for two omniscience
principles:\\

\noindent \textbf{Limited Principle of Omniscience (LPO).} \textit{For
any real number $\alpha\ge0$, either $\alpha=0$ or $\alpha>0$.}\\

\noindent \textbf{Lesser Limited Principle of Omniscience (LLPO).}
\textit{For any real number $\alpha$, either $\alpha\leq0$ or $\alpha\geq0.$}\\

\noindent A statement is considered \textit{constructively invalid}
if it implies an omniscience principle. Following Bishop, we may at
times use the italicized \textit{not} to indicate a constructively
invalid statement.\footnote{For more information concerning Brouwerian counterexamples, and other
omniscience principles, see {[}M83{]} and {[}M89{]}.}

\textit{Normed spaces.} In a real normed linear space $X$, we denote
by $B_{r}(c)$ the closed ball with center $c$ and radius $r\geq0$.
We write simply $B$ for the unit ball $B_{1}(0)$, and $\partial B$
for its boundary. The convex hull of points $u,v\in X$ will be denoted
$\mathscr{H}\{u,v\}$.

\section{Uniform convexity}

\label{sec:unif}

A real normed linear space $X$ is \textit{uniformly convex} if for
every $\varepsilon>0$ there exists $q<1$ such that, if $u,v\in\partial B$
with $\|u-v\|\geq\varepsilon,$ then $\|\frac{1}{2}(u+v)\|\leq q$.

Lemma 1 in {[}R07{]} is stated below as Lemma \ref{Lm1}. This lemma
extends the defining property of uniform convexity, stated for unit
vectors, to arbitrary vectors; it is required for Richman's characterization
of nearly convex subsets of uniformly convex spaces.\footnote{{[}R07{]}, Theorem 4: \textit{A subset of a uniformly convex space
is nearly convex if and only if its closure is convex.}} The proof here fills a minor gap in the proof found in {[}R07{]},
where the unit vectors at the end need not attain the required distance.

\begin{lm} \label{Lm1} (Richman) Let $X$ be a real normed linear
space. If $X$ is uniformly convex, then for any $\varepsilon>0$
there exists $q<1$ with the following property: For any point $c\in X$
and any real number $r\ge0$, if $u,v\in B_{r}(c)$ with $\|u-v\|\geq r\varepsilon$,
then $\|\frac{1}{2}(u+v)-c\|\leq qr$. \end{lm} \textbf{Proof.}\footnote{The main ideas in this proof are drawn from {[}R07{]}.}
Given $\varepsilon>0$, choose $q'<1$ using $\varepsilon'=\varepsilon/2$
in the definition of uniform convexity. Set $\delta=\min\{\varepsilon/12,(1-q')/4,1/3\}$,
and take $q=1-\delta$.

(a) First consider the situation in which $r>0$; we may assume that
$c=0$ and $r=1$. Let $u,v\in B$ with $\|u-v\|\geq\varepsilon$.
Either $\|u\|<1-2\delta$ or $\|u\|>1-3\delta$. In the first case,
$\|\frac{1}{2}(u+v)\|<\frac{1}{2}((1-2\delta)+1)=q$. Similarly, there
are two cases for $v$.

In the fourth case, both norms are $>1-3\delta$. Set $\bar{u}=u/\|u\|,$
$\bar{v}=v/\|v\|,$ $m=\frac{1}{2}(u+v),$ and $m'=\frac{1}{2}(\bar{u}+\bar{v})$.
Thus $\|u-\bar{u}\|=1-\|u\|<3\delta$, and similarly for $v$. Now
\begin{tabbing} \hspace{70pt} $\|\bar{u}-\bar{v}\|$ \= $\geq\|u-v\|-\|u-\bar{u}\|-\|\bar{v}-v\|$
\\
 \> $>\varepsilon-3\delta-3\delta\geq\varepsilon/2=\varepsilon'$
\end{tabbing} so $\|m'\|\leq q'$, and it follows that \begin{tabbing}
\hspace{70pt} $\|m\|$ \= $\leq\|m'\|+\|m-m'\|$ \\
 \> $\leq\|m'\|+\frac{1}{2}(\|u-\bar{u}\|+\|v-\bar{v}\|)$ \\
 \> $<q'+3\delta\leq1-\delta=q$ \end{tabbing} This verifies the
required property when $r>0$.

(b) Now consider the general situation, where $r\ge0$; we may assume
that $c=0$. Let $u,v\in B_{r}(0)$ with $\|u-v\|\geq r\varepsilon$,
and let $e>0$. Either $r>0$, or $r<e$. In the first case, (a) applies,
so $\|\frac{1}{2}(u+v)\|\le qr<qr+e$. In the second case, $\|\frac{1}{2}(u+v)\|\le r<e\le qr+e$.
Thus $\|\frac{1}{2}(u+v)\|\le qr$. $\blacksquare$ \\
 \\
 Also required for {[}R07{]}, Theorem 4, is {[}R07{]}, Lemma 3, stated
below as Lemma \ref{Lm3}. The proof here avoids an incorrect assumption
at the beginning of the proof found in {[}R07{]}.

\begin{lm} \label{Lm3} (Richman) Let $X$ be a uniformly convex
real normed linear space, and let $c$ and $d$ be points of $X$.
Then $\mathrm{diam}(B_{r}(c)\cap B_{s})\rightarrow0$ as $r+s\rightarrow\|c-d\|$.\footnote{The diameter of the intersection may not be amenable to constructive
calculation; the expression \textit{diameter tends to zero} means
only that one may find arbitrarily small upper bounds.} \end{lm} \textbf{Proof.}\footnote{The main ideas in this proof are drawn from {[}R07{]}. The proof here
avoids the incorrect assumption $r+s\leq\|x-y\|+\inf\{r,s\}$ found
in {[}R07{]}. The proof here also removes the restriction to distinct
centers; this would allow a slight simplification in the proof of
{[}R07{]}, Theorem 4.} Set $\rho=\|c-d\|$. Given $\varepsilon>0$, either $\rho<\varepsilon/4$
or $\rho>0$. In the first case, set $\delta=\varepsilon/4$. If $|r+s-\rho|<\delta$,
then $r\leq r+s\leq\rho+\delta<\varepsilon/2$, so $\mathrm{diam}B_{r}(c)\leq2r<\varepsilon$.

In the second case, choose $q<1$ in Lemma \ref{Lm1} using $\varepsilon'=\varepsilon/\rho$,
and set $\delta=\mathrm{min}\{\varepsilon/4,\rho(1-q)\}$. Let $|r+s-\rho|<\delta$,
and let $u,v\in B_{r}(c)\cap B_{s}(d)$. Either $r>\rho-\delta$ or
$r<\rho$. In the first of these subcases, $s<\rho-r+\delta<2\delta$,
so $\|u-v\|\leq2s<4\delta\leq\varepsilon$. Similarly, there are two
subcases for $s$.

There remains only the fourth subcase, where both $r<\rho$ and $s<\rho$.
Suppose $\|u-v\|>\varepsilon$. Then $\|u-v\|>\varepsilon'\rho>\varepsilon'r,$
with $r>0$, and similarly for $s$. Thus $\|\frac{1}{2}(u+v)-c\|\leq qr$
and $\|\frac{1}{2}(u+v)-d\|\leq qs$, so 
\begin{center}
$\rho=\|c-d\|\le q(r+s)<r+s$ 
\par\end{center}

\noindent \begin{flushleft}
and \begin{tabbing} \hspace{70pt} $r+s-\rho$ \= $\ge r+s-q(r+s)$
\\
 \> $=(1-q)(r+s)>(1-q)\rho\ge\delta$ \end{tabbing} which contradicts
the choice of $r$ and $s$; thus $\|u-v\|\leq\varepsilon$. $\blacksquare$
\par\end{flushleft}

\section{Near and metric convexity}

\label{sec:met} A subset $S$ of a metric space $(M,\rho)$ is \textit{nearly
convex} if for any $x,y\in S$, $\lambda>0$, and $\mu>0$, with $\rho(x,y)<\lambda+\mu$,
there exists $z\in S$ such that $\rho(x,z)<\lambda$ and $\rho(z,y)<\mu$.
This version of convexity was introduced in {[}M83{]} in connection
with continuity problems.

Richman {[}R07{]} has defined a subset $S$ of a metric space to be
\textit{metrically convex} if for any $x,y\in S$, $\lambda\geq0$,
and $\mu\geq0$, with $\rho(x,y)=\lambda+\mu$, there exists $z\in S$
such that $\rho(x,z)=\lambda$ and $\rho(z,y)=\mu$. It will follow
from Example \ref{ex-com-pt} and Theorem \ref{com-met} that the
statement ``Every real normed linear space is metrically convex"
is constructively invalid.

The following is a comment in {[}R07{]}; a proof is included here.

\begin{prop} \label{met-near} (Richman) A metrically convex subset
of a metric space is nearly convex. \end{prop} \textbf{Proof.} Let
$x,y\in S$, $\lambda>0$, and $\mu>0$, with $\rho(x,y)<\lambda+\mu$.
Set $\gamma=\lambda+\mu-\rho(x,y)$. Either $\lambda<\gamma$ or $\lambda>\gamma/2$.
In the first case, set $\lambda'=0$ and $\mu'=\mu+\lambda-\gamma$.
Then $0\leq\lambda'<\lambda$, $0\leq\mu'<\mu$, and $\rho(x,y)=\lambda'+\mu'$,
so $z$ may be selected in $S$ using metric convexity. Similarly,
there are two cases for $\mu$.

In the fourth case, we have $\lambda>\gamma/2$ and $\mu>\gamma/2$.
Set $\lambda'=\lambda-\gamma/2$ and $\mu'=\mu-\gamma/2$. Then $0\leq\lambda'<\lambda$,
$0\leq\mu'<\mu$, and $\rho(x,y)=\lambda'+\mu'$, so again metric
convexity provides a suitable element $z\in S$. $\blacksquare$ \\
 \\
 The statement in {[}R07{]}, ``Convex subsets of a normed space are
metrically convex", is incorrect. Convex subsets that are \textit{not}
metrically convex are found nearly everywhere, as is shown by the
following Brouwerian counterexample.

\begin{ex} \label{ex-conv-met} Every nontrivial convex subset $S$
of a real normed linear space $X$ contains a convex subset $T$ such
that the statement ``$\:T$ is metrically convex" is constructively
invalid; the statement implies LLPO. \end{ex} \textbf{Proof.} By
\textit{nontrivial} we mean that $S$ contains at least two distinct
points $u$ and $v$. Given any real number $\alpha$, let $T$ be
the convex hull of the points $\frac{1}{2}(u+v)\pm|\alpha|(u-v)/\|u-v\|$.
Then $T\subset S$ (since we may assume that $|\alpha|$ is small),
and $T$ is isometric with the subset $V=\mathscr{H}\{-|\alpha|,|\alpha|\}$
of $\mathbb{R}$.

Set $x=-|\alpha|$, $y=|\alpha|$, $\lambda=|\alpha|+\alpha$, and
$\mu=|\alpha|-\alpha$. Then $x,y\in V$, $\lambda\geq0$, $\mu\geq0$,
and $\lambda+\mu=\rho(x,y)$. Under the hypothesis ``$\:T$ is metrically
convex", so also is $V$; thus there exists $z\in V$ such that $\rho(x,z)=\lambda$
and $\rho(z,y)=\mu$. Since $-|\alpha|\leq z\leq|\alpha|$, we have
$z+|\alpha|=\rho(x,z)=\lambda=|\alpha|+\alpha$, so $z=\alpha$.

This shows that $\alpha\in\mathscr{H}\{-|\alpha|,|\alpha|\}$, so
there exists $t$ with $0\leq t\leq1$ and $\alpha=(1-t)(-|\alpha|)+t|\alpha|=(2t-1)|\alpha|$.
Note that if $\alpha>0$, then $t=1$, while if $\alpha<0$, then
$t=0$. Either $t<1$ or $t>0$, and it follows that either $\alpha>0$
is contradictory or $\alpha<0$ is contradictory. Thus either $\alpha\le0$
or $\alpha\ge0$, and LLPO results. $\blacksquare$ \\
 \\
 \textit{Convex hulls and intervals.} It is easily seen that the convex
hull $\mathscr{H}\{a,b\}$ is a dense subset of the closed interval
$[a,b]$ whenever $a\le b$, and coincides with the interval when
$a<b$. Peter Schuster has raised the question of whether these sets
are always equal.\footnote{{[}S03{]}, page 446.} A negative answer
was given by Example 10.11 in {[}M83{]}. However, a simpler example
is included in the proof of Example \ref{ex-conv-met} above; although
the real number $\alpha$ lies in the closed interval $[-|\alpha|,|\alpha|]$,
the statement ``$\,\alpha\in\mathscr{H}\{-|\alpha|,|\alpha|\}$ for
all $\alpha\in\mathbb{R}$" implies LLPO. \\
 \\
 Example \ref{ex-conv-met} notwithstanding, we have the following
result for \textit{complete} convex subsets:

\begin{thm} \label{complete-met} A complete convex subset $S$ of
a real normed linear space $X$ is metrically convex. \end{thm} \textbf{Proof.}\footnote{The original proof of Theorem \ref{complete-met} constructed a Cauchy
sequence of points converging to the required point $z$; this method
is more typical of constructive practice. However, constructing such
a sequence requires an appeal to the \textit{axiom of countable choice}.
The nested-sequence method of proof used here, suggested by the referee,
avoids this axiom. Avoiding use of the \textit{axiom of countable
choice} is considered a desirable goal by many constructivists. For
more details concerning this issue, see {[}R01{]} and {[}S03{]}.} Let $x,y\in S$, $\lambda\geq0$, $\mu\geq0$, with $\|x-y\|=\lambda+\mu$;
set $\rho=\lambda+\mu$. When $\rho>0$, we write $w=(\mu x+\lambda y)/\rho\,$;
in this situation we have $\|w-x\|=\lambda$ and $\|w-y\|=\mu$. For
all positive integers $n$, define 
\begin{center}
$S_{n}=\{x:\rho<1/n\}\cup\{w:\rho>0\}$ 
\par\end{center}

Let $m<m'$, and $y\in S_{m'}$. If $y=x$, then $\rho<1/m'<1/m$,
so $y\in S_{m}$, while if $y=w$, with $\rho>0$, then again $y\in S_{m}$;
thus $S_{m'}\subset S_{m}$. Also, if $x,w\in S_{n}$, with $\rho>0$,
then $\|w-x\|=\lambda\le\rho<1/n$; thus $\mathrm{diam}\;S_{n}<1/n$.

Thus the sequence $\{\overline{S_{n}}\}_{n}$ of closures is a nested
sequence of nonvoid closed subsets of $S$ with diameters tending
to zero. Since $S$ is complete, there is a unique point $z$ in the
intersection of this sequence.

To show that the point $z$ satisfies the required conditions, let
$\varepsilon>0$. If $\rho>0$, then $S_{n}=\{w\}$ eventually, so
$z=w$ and $\|z-x\|=\lambda$. On the other hand, if $\rho<\varepsilon$,
then $\lambda<\varepsilon$ , and $\|s-x\|<\varepsilon$ for all $s\in S_{1}$,
so $\|z-x\|\le\varepsilon$. It follows that $|\,\|z-x\|-\lambda\,|\le\varepsilon$.
Thus $\|z-x\|=\lambda$. Similarly, $\|z-y\|=\mu$. Thus $S$ is metrically
convex. $\blacksquare$

\begin{cor} \label{x-met} Every complete real normed linear space
$X$ is metrically convex. \end{cor}

\noindent In the proof of Theorem \ref{complete-met}, the passage
from the sequence $\{S_{n}\}_{n}$ to the sequence $\{\overline{S_{n}}\}_{n}$
of closures is necessary. Although each set $S_{n}$ contains at most
two elements, these sets are \textit{not} closed. To demonstrate this,
consider any real number $\alpha\ge0$, and take $X=\mathbb{R},\;x=0,\;y=\alpha,\;\lambda=\alpha$,
and $\mu=0$. Then each set $S_{n}$ is included in the following
Brouwerian counterexample.

\begin{ex} The statement ``In a complete real normed linear space,
any nonvoid subset containing at most two elements is closed" is
constructively invalid; the statement implies LPO. \end{ex} \textbf{Proof.}
Given $\alpha\ge0$, choose any real number $c>0$, and define 
\begin{center}
$S=\{0:\alpha<c\}\cup\{\alpha:\alpha>0\}$ 
\par\end{center}

To show that $\alpha$ lies in the closure of $S$, let $\varepsilon>0$;
we may assume that $\varepsilon<c$. Either $\alpha<\varepsilon$
or $\alpha>0$. In the first case, 0 is a point of $S$ with distance
to $\alpha$ less than $\varepsilon$. In the second case, $\alpha\in S$,
and this suffices. Thus $\alpha\in\overline{S}$. By hypothesis, $\alpha\in S$,
so $\alpha$ must lie in one of the sets forming the union. Hence
LPO results. $\blacksquare$

\section{Osculation}

\label{sec:osc}

Two closed balls, $B_{r_{1}}(c_{1})$ and $B_{r_{2}}(c_{2})$, in
a real normed linear space $X$ are said to be \textit{osculating}
if $\|c_{1}-c_{2}\|=r_{1}+r_{2}$. When $r_{1}+r_{2}>0$, the balls
are \textit{nondegenerate}.

Theorem 5 in {[}R07{]} considers various conditions on $X$ which
ensure that a point common to two osculating balls is \textit{unique}.
On the other hand, while the \textit{existence} of a common point
is evident in the case of nondegenerate osculating balls, a common
point does \textit{not} always exist, as the Brouwerian counterexample
below will show. Thus the implication ``(2) implies (1)" in {[}R07{]},
Theorem 5, must be understood to involve only uniqueness.

\begin{ex} \label{ex-com-pt} The statement ``Osculating balls in
a real normed linear space $X$ always have at least one common point"
is constructively invalid; it implies LLPO. \end{ex} \textbf{Proof.}
Given any real number $\alpha$, set $X=\mathbb{R}\alpha$, the linear
subspace of $\mathbb{R}$ generated by the single point $\alpha$,
with the induced norm.\footnote{For an extensive discussion of this space, see {[}R82{]}.}
Set $\alpha^{+}=\max\{\alpha,0\}$ and $\alpha^{-}=\max\{-\alpha,0\}$.
Since $\alpha^{+}+\alpha^{-}=|\alpha|$, the balls $B_{\alpha^{-}}(0)$
and $B_{\alpha^{+}}(\alpha)$ in $X$ are osculating.

By hypothesis, there exists a point $x$ of $X$ in the intersection
$I$ of these balls, with $x=\lambda\alpha$ for some $\lambda\in\mathbb{R}$.
Note that if $\alpha>0$, then $I=B_{0}(0)\cap B_{|\alpha|}(\alpha)$,
so $x=0$ and $\lambda=0$, while if $\alpha<0$, then $I=B_{|\alpha|}(0)\cap B_{0}(\alpha)$,
so $x=\alpha$ and $\lambda=1$. Either $\lambda>0$ or $\lambda<1$,
and it follows that either $\alpha>0$ is contradictory or $\alpha<0$
is contradictory. Thus either $\alpha\le0$ or $\alpha\ge0$, and
LLPO results. $\blacksquare$ \\
 \\
 It will be convenient to state for reference a few basic properties
of osculating balls, all of which are included in {[}R07{]}.

\begin{lm} \label{lm-osc} Let $B_{r_{1}}(c_{1})$ and $B_{r_{2}}(c_{2})$
be osculating balls in a real normed linear space.

(a) If $v$ is a common point, then $v$ lies on the boundary of each
ball.

(b) If $u$ and $v$ are common points, then all points in the convex
hull $\mathscr{H}\{u,v\}$ are also common points.

(c) When $r_{1}+r_{2}>0$, the point $z=(r_{2}c_{1}+r_{1}c_{2})/(r_{1}+r_{2})$
is a common point.

(d) Let $r_{1}+r_{2}>0$. If a common point $v$ lies on the line
through the centers, then $v$ is the point $z$ specified in (c).
\end{lm} \textbf{Proof.} In (a), we have $\|v-c_{i}\|\leq r_{i}$.
Also, 
\begin{center}
$\|v-c_{1}\|\ge\|c_{1}-c_{2}\|-\|v-c_{2}\|\ge(r_{1}+r_{2})-r_{2}=r_{1}$ 
\par\end{center}

Thus $\|v-c_{1}\|=r_{1}$, and $v\in\partial B_{r_{1}}(c_{1})$. Similarly,
$v\in\partial B_{r_{2}}(c_{2})$. The remaining properties follow
from simple calculations. $\blacksquare$ \\
 \\
 From the definition of metric convexity and Lemma \ref{lm-osc}(a)
we have the following.

\begin{thm} \label{com-met} Let $X$ be a real normed linear space.
Osculating balls in $X$ always have at least one common point if
and only if $X$ is metrically convex. \end{thm} From this theorem
and Corollary \ref{x-met} we obtain the following.

\begin{cor} \label{complete-com} In a complete real normed linear
space, osculating balls always have at least one common point. \end{cor}

\section{Strict convexity}

A real normed linear space $X$ is \textit{strictly convex} if the
convex hull $\mathscr{H}\{u,v\}$ of two points is contained in the
boundary $\partial B$ of the unit ball $B$ only when $u=v$. It
follows that for any ball $B_{r}(c)$, if $\mathscr{H}\{u,v\}\subset\partial B_{r}(c)$,
then $u=v$.

The following is Richman's theorem relating strict convexity to the
uniqueness of points of osculation.

\begin{thm} (Richman)\footnote{{[}R07{]}, Theorem 5. Only the first five conditions are listed here;
the sixth condition concerns a problem that remains open.} For any real normed linear space $X$, the following conditions are
equivalent.

(1) Any two osculating balls have at most one common point.

(2) Any two nondegenerate osculating balls have at most one common
point.

(3) Any two osculating unit balls have at most one common point.

(4) $X$ is strictly convex.

(5) If $x$ and $y$ are points of $X$ with $\|x+y\|=\|x\|+\|y\|\neq0$,
then $x$ and $y$ are linearly dependent. \end{thm} \textbf{Proof.}\footnote{The proof here is essentially the same as in {[}R07{]}, with several
simplifications.} That (1) implies (2) implies (3) is self-evident.

(3) implies (4). Let $\mathscr{H}\{u,v\}\subset\partial B$. Then
$\|\frac{1}{2}(u+v)\|=1$, so $B_{1}(0)$ and $B_{1}(u+v)$ are osculating
unit balls. Since $u$ and $v$ both lie in each ball, $u=v$.

(4) implies (1). Let $B_{r_{1}}(c_{1})$ and $B_{r_{2}}(c_{2})$ be
osculating balls, and let $u$ and $v$ be common points. By Lemma
\ref{lm-osc}(a)(b), $\mathscr{H}\{u,v\}\subset\partial B_{r_{1}}(c_{1})$,
so $u=v$.

(1) implies (5). Given $\|x+y\|=\|x\|+\|y\|>0$, the balls $B_{\|x\|}(x)$
and $B_{\|y\|}(-y)$ are osculating and nondegenerate, so by Lemma
\ref{lm-osc}(c) they have a common point of the form $z=\lambda x+(1-\lambda)(-y)$.
Since 0 is clearly a common point, $z=0$, and this results in a linear
dependence relation for $x$ and $y$.

(5) implies (3). Let $B_{1}(c_{1})$ and $B_{1}(c_{2})$ be osculating
unit balls, and let $v$ be any common point. Then 
\begin{center}
$\|(c_{1}-v)+(v-c_{2})\|=\|c_{1}-c_{2}\|=2=\|c_{1}-v\|+\|v-c_{2}\|$ 
\par\end{center}

\noindent so $c_{1}-v$ and $v-c_{2}$ are linearly dependent. We
have then $\lambda(v-c_{1})+\mu(v-c_{2})=0$, where either $\lambda\neq0$
or $\mu\neq0$. We may assume that $\mu\neq0$, and furthermore that
$\mu=1$; thus $(\lambda+1)v=\lambda c_{1}+c_{2}$. It follows that
$(\lambda+1)(v-c_{1})=c_{2}-c_{1}$; taking norms here, we find that
$|\lambda+1|=2$, so $\lambda$ is either -3 or 1. Also, $(\lambda+1)(v-c_{2})=\lambda(c_{1}-c_{2})$;
taking norms here, we have $|\lambda|=1$, so -3 is ruled out, and
it follows that $v=\frac{1}{2}(c_{1}+c_{2})$. Thus $v$ is unique.
$\blacksquare$\\
\newpage{}

\section*{{\large{}References}}

\noindent {\small{}{[}B67{]} Errett Bishop, }\textit{\small{}Foundations
of Constructive Analysis.}{\small{} McGraw-Hill Book Co., New York-Toronto-London,
1967.}{\small\par}

\noindent {\small{}{[}M83{]} Mark Mandelkern, }\textit{\small{}Constructive
Continuity.}{\small{} Memoirs Amer. Math. Soc. 277 (1983).}{\small\par}

\noindent {\small{}{[}M85{]} \rule[2pt]{20pt}{0.5pt}, }\textit{\small{}Constructive
mathematics,}{\small{} Math. Mag. 58 (1985), 272-280.}{\small\par}

\noindent {\small{}{[}M89{]} \rule[2pt]{20pt}{0.5pt}, }\textit{\small{}Brouwerian
counterexamples,}{\small{} Math. Mag. 62 (1989), 3-27.}{\small\par}

\noindent {\small{}{[}R82{]} Fred Richman, }\textit{\small{}Meaning
and information in constructive mathematics}{\small{}, Amer. Math.
Monthly 89 (1982), 385-388.}{\small\par}

\noindent {\small{}{[}R99{]} \rule[2pt]{20pt}{0.5pt}, }\textit{\small{}Existence
proofs}{\small{}, Amer. Math. Monthly 106 (1999), 303-308.}{\small\par}

\noindent {\small{}{[}R01{]} \rule[2pt]{20pt}{0.5pt}, }\textit{\small{}Constructive
mathematics without choice,}{\small{} Reuniting the antipodes - constructive
and nonstandard views of the continuum (Venice, 1999), 199-205, Synthese
Lib., 306, Kluwer Acad. Publ., Dordrecht, 2001.}{\small\par}

\noindent {\small{}{[}R07{]} \rule[2pt]{20pt}{0.5pt}, }\textit{\small{}Near
convexity, metric convexity and convexity}{\small{}, Rocky Mountain
J. Math. 37 (2007), 1305-1314.}{\small\par}

\noindent {\small{}{[}S03{]} Peter M. Schuster, }\textit{\small{}Unique
existence, approximate solutions, and countable choice,}{\small{}
Theoretical Computer Science 305 (2003), 433-455. }{\small\par}

{\small{}\vspace{1cm}
}{\small\par}
\noindent \begin{flushleft}
{\small{}New Mexico State University }\\
{\small{}Las Cruces, New Mexico, USA }\\
\textit{\small{}e-mail:}{\small{} mmandelk@nmsu.edu}\\
\textit{\small{}e-mail:}{\small{} mandelkern@zianet.com}\\
\textit{\small{}web:}{\small{} www.zianet.com/mandelkern}\\
{\small{} July 31, 2010 }{\small\par}
\par\end{flushleft}
\end{document}